\documentclass[11pt]{article}

\UseRawInputEncoding 

\newtheorem{theorem}{Theorem}[section]
\newtheorem{lemma}[theorem]{Lemma}

\newtheorem{definition}[theorem]{Definition}

\newcommand\qed{\begin{flushright} {\bf q.e.d.} \end{flushright} }
\newcommand\prf{\noindent {\bf Proof :}}  

\newcommand\bits{\{0,1\}}
\newcommand\uu{{\bits^*}}

\newcommand\bk{{\bf k}}
\newcommand\bs{{\bf s}}
\newcommand\bl{{\bf {\ell}}}

\newcommand\fla{{\sf Fla}}
\newcommand\defc{{\sf Def}}
\newcommand\vba{\mbox{BA}\ \vdash}
\newcommand\simba{\sim_{\mbox{BA}}}
\newcommand\eqcd{{\mbox{Eq}_{C,D}}}

\begin{document}

\title{Extended Frege proofs, circuits and rewriting} 

\author{Jan Kraj\'{\i}\v{c}ek}

\date{Faculty of Mathematics and Physics\\
Charles University\thanks{Sokolovsk\' a 83, Prague, 186 75,
The Czech Republic, {\tt jan.krajicek@protonmail.com}}}

\maketitle

\begin{abstract}
Inspired by a statement about Extended Frege proof systems by
Jain and Jin \cite[Lemma 3]{JJ22} we prove that:
\begin{itemize} 

\item there is a p-time binary relation $\approx$ between circuits
that implies their logical equivalence,

\item if the equivalence $C \equiv D$ 
has an Extended Frege or a Circuit Frege proof 
with $k$ steps  
then there is a chain of circuits $E_i$
$$
C = E_0 \approx \dots \approx E_t = D
$$
with $t \le O(k \log k)$, and the total size of the chain
is bounded above by $(k+|C|+|D|)^{O(1)}$.
\end{itemize}
\end{abstract}


\section*{Introduction}

Jain and Jin \cite{JJ22} 
constructed a new type of cryptographic
indistinguishability obfuscator for pairs of circuits 
whose functional equality has short Extended Frege (EF) proofs. 
For this they established a certain property of EF proofs that 
can be informally described as follows.
There is a p-time algorithm that from a circuit
$C$ and $1^{(s)}$ computes another circuit $C'$ (a padded version of $C$)
such that $C \equiv C'$ and whenever the equivalence 
$C \equiv D$ has a size $s$ EF-proof then $C'$ can be converted into
$D'$ by a polynomially long chain of intermediate circuits. 
All these intermediate circuits have the same 
underlying dag and are obtained from $C'$ in consecutive steps,
each step changing the functionality of a logarithmic number of gates 
scattered in the circuit
in a way that the subcircuits they induce (as defined there) compute the same 
Boolean function. 

Here we shall use one of the ideas behind 
the finite Boolean valuation characterization of EF proof size 
given in \cite{Kra-frege} to prove a statement akin to \cite[L.3]{JJ22}. 
That idea is to consider the number of steps in Frege proofs instead of EF proof size 
and to use equational logic to analyze Frege proofs with a bounded number of steps.
The resulting Theorem \ref{main}, linking proofs of $C \equiv D$
with a rewriting of $C$ into $D$, isolates a rather rudimentary
property of Extended Frege and Circuit Frege proofs.

The paper is organized as follows.
After some preliminaries about proof complexity in Sec.\ref{prelim}
we move the construction to equational logic and rewriting 
in Sec.\ref{eqlo}. The theorem is proved in 
Sec.\ref{re} and we conclude by a few remarks there.
The reader requiring a more detailed proof complexity background 
can find it in \cite{prf}.

\section{Proof complexity} \label{prelim}

We shall consider formulas and circuits using the DeMorgan language: 
$0,1, \neg$ and binary $\vee, \wedge$. Formulas are terms in this language.
A size $s$ 
circuit $C$ with inputs $x_1,\dots, x_n$ is specified by a straight line program 
giving instructions on how to compute values of gates $y_j$, $j = 1, \dots, s$, from values
of earlier gates: 
$y_j$ is defined to be one of the inputs or a constant, or is computed from
$j_1 < j$ or from $j_1 < j_2 < j$ by $\neg$ or by one of $\vee,\wedge$, respectively.
The program determines a unique directed acyclic graph (dag) whose edges (wires) into a gate
go from  the gates used in its definition. In particular, gates defined as inputs or constants
have in-degree $0$ and all other have in-degree $1$ or $2$. 
We assume that the last gate $y_s$ is the
designated output gate of $C$ and that it is the only out-degree $0$ gate.
Note that the program determines also an ordering of the 
wires so we can talk about first (left)
or second (right) wire entering a gate.

A {\em Frege proof system} in any complete propositional language is given by a finite, 
sound and implicationally complete set of Frege rules. For definiteness we
shall fix any system F in the DeMorgan language that has a finite
number of axiom schemes and modus ponens
$$
	\frac{\alpha_1\ ,\ \neg \alpha_1 \vee \alpha_0}{\alpha_0}
$$
as the only inference rule. This choice simplifies a bit the argument for 
Claim 2 in the proof of Theorem \ref{main}
but it is without a loss of generality as 
by Reckhow's theorem all Frege system p-simulate each other, cf.
\cite{CooRec} or \cite[Sec.2.3]{prf}.

We shall need the following basic result about Frege systems. 
A proof is {\em tree-like} if every formula is used as a 
hypothesis of at most one inference.

\begin{theorem} [\cite{Kra-lower}] \label{treethm}
{\ }

If a formula has an F-proof with $k$ steps then it has a tree-like
F-proof with $O(k \log k)$ steps.
\end{theorem}
The idea behind its proof is to transform an original proof with steps $\alpha_i$
step-by-step into a proof of formulas 
$\alpha_i^* := \bigwedge_{j \le i}\alpha_j$:
$\alpha^*_i$ contains all hypotheses needed to derive $\alpha_{i+1}$ 
and hence also $\alpha^*_{i+1}$, and thus the resulting proof is tree-like.
Moreover one can combine the proofs of the implications
$\alpha^*_i \rightarrow \alpha^*_{i+1}$ in a binary tree fashion so that the constructed
proof tree is balanced. See \cite[Thm.2.2.1]{prf} for details of the proof.

\bigskip

System F operates with formulas but any 
proof system can express statements about circuits using their definition. 
Let
$\defc_C(x_1, \dots, x_n, y_1, \dots, y_s)$ be the conjunction of all instructions
defining $C$; it can be written as a CNF. If $D(\overline x)$ 
is a size $t$ circuit then we can express the functional equivalence of $C$ and $D$
by a formula to be denoted $\eqcd$:
\begin{equation}\label{eq}
[\defc_C(\overline x, \overline y) \wedge \defc_D(\overline x, \overline z)]\ 
\rightarrow\ y_s \equiv z_t\ .
\end{equation}
It is another matter to be able to use circuits as lines in a proof. For this
Cook and Reckhow \cite{CooRec} strengthened F to the {\em Extended Frege system} 
EF by augmenting F 
with Tseitin's \cite{Tse68} {\em extension rule}. This rule 
allows to introduce  during a proof a new {\em extension axiom}
$$
r \equiv \beta
$$
if the {\em extension variable} $r$ does not occur earlier in the proof and does not 
occur in $\beta$ or in the last formula
being proved. Using the extension axioms one may define a circuit and thus EF can
indirectly operate with circuits.

However, it is more elegant to allow circuits directly
as lines. This was formally defined by Je\v r\' abek \cite{Jer04} 
as {\em Circuit Frege system} CF.
For a circuit $C$ let $\fla(C)$ be the unique formula obtained by unwinding $C$:
$\fla(C) := C$ if $C$ is an input variable of a constant, $\fla(\neg C):=\neg \fla(C)$
and $\fla(C \circ D):= \fla(C)\circ \fla(D)$ for $\circ = \vee, \wedge$ (the straight line program determines the order of wires and so also the order of $C, D$ in $C\circ D$).
Note that the size of $\fla(C)$ can be exponentially bigger that the size $|C|$
of $C$
but the number of different subformulas occurring in $\fla(C)$ is bounded above 
$|C|$ (in fact, by the number of different subcircuits in $C$). 

CF is defined as F allowing circuits and not just formulas as its lines
and it has one new inference rule:
$$
\frac{C}{D}\ ,\ \mbox{assuming } \ \fla(C) = \fla(D)\ .
$$
The property of $C,D$ that $\fla(C) = \fla(D)$ is p-time decidable
(cf. \cite{Jer04}) 
and CF is thus a proof system in the sense of \cite{CooRec}. 
Je\v r\' abek \cite{Jer04} proved that
EF and CF polynomially simulate each other w.r.t. proofs of formulas
(cf. also \cite[L.7.2.2]{prf}).

For $P$ one of EF or CF and $\tau$ a tautology define $\bs_P(\tau)$ and $\bk_P(\tau)$
to be the minimal size of a $P$-proof of $\tau$ and the minimal number of steps
in a $P$-proof of $\tau$, respectively. 
Using the relations between $\bk_{EF}, \bk_{F}$ and $\bs_{EF}$ from \cite{CooRec}
(cf. \cite[Sec.2.5]{prf}) it holds that the quantities:
\begin{equation}\label{eq2}
\bk_{EF}(\eqcd) + |C| + |D|\ ,\ 
\bk_{F}(\eqcd) + |C| + |D|\ ,\ 
\end{equation}
\begin{equation} \label{eq3}
\bs_{EF}(\eqcd)\ ,\ 
\bs_{CF}(\eqcd)\ ,\ 
\bs_{CF}(C \equiv D)
\end{equation}
are polynomially related.
Here the equivalence $C \equiv D$
is an abbreviation for the circuit $(C\wedge D) \vee (\neg C \wedge \neg D)$.

\medskip

For a formula $A$ and an $F$-proof $\pi$ let $\bl(A)$ and $\bl(\pi)$ denote 
{\em the number of distinct formulas} that occur
as subformulas in $A$ and in $\pi$, respectively,
and let $\bl_F(A)$ by the minimum 
of $\bl(\pi)$ over all $F$-proofs
of $\pi$ of $A$. Obviously:
$$
\bk_F(A) \le \bl_F(A) \le \bs_F(A)\ .
$$
In fact, a stronger bound holds.

\begin{lemma} [{\cite[L.4.4.6]{kniha}}] \label{8.10.26a}
	{\ }
	
If a formula $A$ has an (tree-like) F-proof with $k^*$ steps then it has also
an (tree-like) F-proof $\pi^*$ with $k^*$ steps such that
$$
\bl_F(\pi^*) = O(k^*) + \bl(A)\ .
$$
	
\end{lemma}
The lemma was derived from \cite[L.4.4.4]{kniha} that states
that any $F$-proof $\pi$ is a substitution instance
of a proof $\tilde \pi$ in which every formula has the logical depth bounded above by
$O(\bk(\pi))$. This fact was proved first using a more general result about predicate calculus
from \cite{Kra-steps} and later another proof using a term unification algorithm
was given in \cite{KraPud88}. See also \cite[Sec.1]{Kra-frege} or \cite[Sec.2.5]{prf}.

\section{Equational logic and rewriting} \label{eqlo}

The following set of equations BA between terms in the DeMorgan language axiomatizes 
Boolean algebras:
\begin{itemize}

\item $\neg 0 = 1$, $\neg 1 = 0$, $0\vee a = a$, $1\wedge a = a$,

\item $a\vee \neg a = 1$, $a \wedge \neg a = 0$,

\item $a\vee b = b \vee a$, $a\wedge b = b \wedge a$, 

\item $a \vee (b \vee c) = (a \vee b) \vee c$, 
$a \wedge (b \wedge c) = (a \wedge b) \wedge c$, 

\item $a \vee (b \wedge c) = (a \vee b) \wedge (a\vee c)$, 
 $a \wedge (b \vee c) = (a \wedge b) \vee (a\wedge c)$.

\end{itemize}
For two terms $s,t$ we shall write 
$$
\vba s = t
$$
iff $s=t$ can be derived from BA in 
equational logic having the following inference rules
($u,v,w$ are terms):

\begin{itemize}
\item reflexivity, symmetry and transitivity:
	$$
	\frac{}{u=u}\ ,\ \ \ \ \frac{u=v}{v=u}\ ,\ \ \ \ 
	\frac{u=v\ \ v=w}{u=w} 
	$$
	
\item 	congruences:
$$
\frac{u=v}{\neg u = \neg v}\ ,\ \ \ \ 
\frac{u=v\ \ u'=v'}{u\circ u' = v \circ v'}\ ,
$$
for $\circ = \vee, \wedge$,

\item substitution:
$$
\frac{u = v}{\sigma(u) = \sigma(v)}\ 
$$
where $\sigma$ is a substitution.
\end{itemize}

A derivation of an equation from BA is {\em tree-like} if each equation is used at most
once as a hypothesis of an inference.

\begin{lemma} \label{8.10.26b}
	{\ } 
	
Assume that there is a tree-like $F$-proof $\pi^*$ of formula $A$ with $k^*$ steps
in which at most $\ell^* $ different formulas occur as subformulas. 

Then the equation $A = 1$ has a tree-like proof with $O(k^*)$ steps
in equational logic from axioms BA in which at most 
$(\ell^*)^{O(1)}$ 
different formulas occur as subterms. 

If $A$ is of the form $A_1 \equiv A_2$ then the same is true for the equation $A_1=A_2$.
\end{lemma}

\prf

The equational proof is obtained 
by translating step-by-step $\pi^*$ into a tree-like equational logic  
derivation from BA, each step $B$ in $\pi^*$ being translated into the
equation $B = 1$.

For any axiom $B$ of F there is a constant size tree-like BA-derivation 
of $B = 1$ because all tautologies
get value $1$ in any Boolean algebra, BA axiomatizes them and the equational logic
is complete (Birkhoff's theorem \cite{Bir35}). If we substitute in this fixed derivation 
for variables formulas occurring in $\pi^*$ in the instance of $B$
we create possibly new subformulas 
that do not occur in $\pi^*$. However, there is as most $(\ell^*)^{O(1)}$ of them 
where the constant $O(1)$ depends on the number of occurrences 
of variables in the fixed BA derivation.

To simulate modus ponens we just need
to find a tree-like BA-derivation of $b=1$ from $a=1$ and $\neg a \vee b = 1$
in which each of the two latter formulas occur at most once (otherwise the 
translation of nested applications of modus ponens would grow exponentially). 
An example of such a derivation is: 
\begin{itemize}
	
	\item from the first hypothesis 
	$a=1$ derive $\neg a = \neg 1$ and using the transitivity rule with axiom 
	$\neg 1 = 0$ derive $\neg a = 0$,
	
	\item combine that with the axiom $b = b$ by the congruence rule to get
	$\neg a \vee b = 0 \vee b$ which yields (with axiom $0\vee b = b$) 
	$\neg a \vee b = b$,
	
	\item combine this with the second hypothesis $\neg a \vee b = 1$ to get $b = 1$.
	
\end{itemize}
Analogously as above this produces at most $(\ell^*)^{O(1)}$ new formulas.

\smallskip

For the second part of the lemma 
take any fixed tree-like BA-derivation of $a = b$ from 
$(a\equiv b) = 1$ and append it, after substituting $a : = A_1$ and
$b := A_2$, to the derivation we constructed above.

\qed

Crucially for us the provability of an equation can be characterized by the notion
of a {\em local replacement of equals by equals}. Define the relation
$$
s \simba t
$$
by the condition: there exist terms $w(y), u, v$ where $y$ does not 
occur in $s, t$, and a substitution $\sigma$ for variables in $u, v$ such that
it holds:
\begin{itemize}
	
\item one of the equations $u = v$ or $v = u$ is in BA,	

\item $y$ has one occurrence in $w$,

\item term $s$ is syntactically identical  to $w(y/\sigma(u))$,

\item term $t$ is syntactically identical to $w(y/\sigma(v))$.

\end{itemize}
Note that the relation is 
symmetric and denote by $\simba^*$ its reflexive and transitive closure. 
By {Birkhoff \cite{Bir35}} (cf. also \cite[Thm.3.1.12]{BaaNip}) it holds that
$$
\vba s = t \ \ \ \mbox{ iff }\ \ \ s \simba^* t\ .
$$
The following lemma states additional properties that we shall need.

\begin{lemma} \label{rlemma} 
{\ }	

Assume that the equation $s = t$ has a tree-like derivation from BA with $k'$ steps
in which at most $\ell'$ different terms occur as subterms.
Then there are $c \le k'$, terms $r_0, \dots, r_c$ and  set of term $\Gamma$
of cardinality $|\Gamma| \le (k')^2 + \ell'$ 
such that
$$
s = r_0, \ \dots\ ,\ r_i \simba r_{i+1}\ ,\ \dots\ ,\ 
r_c = t\ 
$$ 	
and such that any term occurring as a subterm 
in the chain $\{r_i\}_i$ is in $\Gamma$.

\end{lemma}

\prf

Fix a tree-like derivation $\pi'$ of $s = t$ 
satisfying the conditions in the hypothesis of the lemma. For an equation $u = v$ in it 
denote by
$\bk(u=v)$ the number of nodes in the tree of the subderivation ending with
$u=v$.  Construct first the chain $\{r_i\}_i$  for all equations $u=v$
in the derivation by induction on $\bk(u=v)$.

As an example of an induction steps assume that
$u_1 \vee u_2 = v_1 \vee v_2$ is derived from 
$u_1=v_1$ and $u_2=v_2$ which were proved by
subproofs with $k_1$ and $k_2$ steps, respectively. By the induction
hypothesis there are chains of terms $\{r^1_i\}$ and $\{r^2_i\}$
of lengths $\le k_1$ and $\le k_2$
linking $u_1$ with $v_1$ and $u_2$ with $v_2$, respectively.

First use $\{r^1_i\}$ to link $u_1 \vee u_2$ with 
$v_1 \vee u_2$ and then link this with $v_1 \vee v_2$ by using the second
chain $\{r^2_i\}$. The combined chain has the length $\le k_1 + k_2 + 1$ and
by the tree-likeness $k_1 + k_2 + 1 = \bk(u_1 \vee u_2 = v_1 \vee v_2)$.

\medskip

To establish the second property let $\Gamma_0$
be the set of terms occurring as subterms in $\pi'$, i.e. $|\Gamma_0|\le\ell'$. 
By induction on $\bk(u=v)$ as before we shall construct set of terms
$\Gamma_{u=v}$ such that all terms occurring in the chain  
$\{r_i\}_i$ linking $u$ and $v$ are among $\Gamma_0 \cup \Gamma_{u=v}$. 

In the case of congruences as above the only new
terms not in $\Gamma_{u_1=v_1}$ and $\Gamma_{u_2=v_2}$ 
are $r^1_i\circ u_2$
(as all subterms of $u_2$ are in $\Gamma_0$) and $v_1 \circ r^2_i$, i.e. at most
$\bk(u=v)$ of them. Add these 
to $\Gamma_{u_1=v_1}\cup \Gamma_{u_2=v_2}$
to form $\Gamma_{u=v}$.

In the case of the substitution rule we can take
$$
\Gamma_{\sigma(u)=\sigma(v)}\ :=\ 
\{\sigma(t)\ |\ t \in  \Gamma_{u=v}\}
$$
as all other subterms of $\sigma(r_i)$ are subterms of substitution instances of variables
and hence also subterms of $\sigma(u), \sigma(v)$ and thus already on $\Gamma_0$.

Altogether we need to add at most $k'$ new terms for each congruence inference
and there are at most $k'$ of them.

\qed

\section{From proofs to chains of circuits} \label{re}

\begin{definition}
	{\ }
	
	Define a binary relation $A \approx B$ between two circuits $A, B$ by: 
	$$
	A \approx B\ \mbox{ iff }\ \fla[A] \simba \fla[B]\ .
	$$	
\end{definition}

\begin{lemma} \label{ptime}
	{\ }
	
	The relation $\approx$ is p-time decidable.
\end{lemma}

\prf

It suffices to show that 
the negated relation $C \not\approx D$ is in the class NL (non-deterministic log-space).
We have that $C \not\approx D$ iff

\begin{itemize}
	
	\item either there is a path $p \in \uu$ such that $C, D$ differ at subcircuits
	defined by $p$ but 
	
	\begin{itemize}
		
		\item not by a pair $u, v$ with the required properties (this requires to check extensions
		of $p$ by at most three steps), or
		
		\item they differ by a suitable pair $u, v$ but there is a path starting at some input gate
		$z_i$ of both $u$ and $v$ that encounters a difference in $C$ and $D$,
		
	\end{itemize}
	
	\item or there are two different paths $p, p'$ such that $C,D$ differ at the corresponding
	positions by suitable pairs $u,v$ and $u', v'$.
	
\end{itemize}
Choosing paths correspond to non-deterministic choices of the next gate
and proceed in parallel in both $C$ and $D$.

\qed

We shall state the bound in the next theorem in terms of $\bk_{CF}$;
as $\bk_{CF}$ is majorized by $\bs_{CF}$ such a formulation implies 
similar statements for all the measures listed in (\ref{eq2}) and (\ref{eq3}).

\begin{theorem} \label{main}
	{\ }
	
Assume that the equivalence $C \equiv D$ has
a CF-proof with $k$ steps. Then there are $t \le O(k \log k)$ and 
circuits $E_i$, $0 \le i \le t$, such that:
\begin{itemize}

\item $E_i \equiv C$, all $i$,

\item $E_0 = C$ and $E_t = D$,

\item $E_i \approx E_{i+1}$ for all $i < t$,

\item the total size of the chain $\{E_i\}_i$ is bounded above
by $(k + |C| + |D|)^{O(1)}$.

\end{itemize}	

\end{theorem}

\prf

Let $\pi: A_1, \dots, A_k$ be a CF-proof of $C \equiv D$ with $k$ steps.

\medskip
\noindent
{\bf Claim 1:} {\em $\pi': \fla[A_1], \dots, \fla[A_k]$ is an $F$-proof of 
$\fla[C]\equiv \fla[D]$.}

\medskip

Theorem \ref{treethm} and Lemma \ref{8.10.26a} imply 

\medskip
\noindent
{\bf Claim 2:}{\em 
There is a tree-like 
F-proof $\pi^*$
of $\fla[C]\equiv \fla[D]$ with \newline
$k^* := O(k \log k)$
steps in which at most 
$\ell^* := O(k^*) + |C| + |D|$ different formulas occur as subformulas\ .
}

\medskip

Lemma \ref{8.10.26b} yields

\smallskip
\noindent
{\bf Claim 3:} {\em  There is a tree-like equational logic derivation of the equation 
$$
\fla[C] = \fla[D]
$$
from BA with $O(k^*)$ steps in which at most $(\ell^*)^{O(1)}$ different formulas 
occur as subterms.}

\smallskip

Claim 3 and Lemma \ref{rlemma} imply

\smallskip
\noindent
{\bf Claim 4:} {\em  There are $t \le O(k^*)$ and formulas 
$\psi_i$, $0 \le i\le t$, such that
$$
\fla[C] = \psi_0\ , \ \dots\ ,\ \psi_i \simba \psi_{i+1}\ ,\ \dots\ ,\ 
\psi_t = \fla[D]
$$ 	
and such that at most 
$$
O((k^*)^2) +  
(\ell^*)^{O(1)} = 
(k + |C| + |D|)^{O(1)}
$$
different subformulas occur as subterms in the chain $\{\psi_i\}_i$.
}

\medskip

The theorem now follows from Claim 4: put $E_0 := C$, $E_t := D$ and for 
$E_i$, $0 < i < t$, take the circuit obtained from $\psi_i$ by identifying 
different occurrences of identical subformulas.

\qed

I do not know if 
the chain with the properties stated in the theorem
can be constructed in polynomial time from the starting CF-proof of the equivalence
$C \equiv D$.
The construction in Jain and Jin \cite{JJ22} is polynomial time but the chain
constructed there does not connect $C$ and $D$ but their padded versions.

\medskip

Let us remark that we
may consider the sequence of circuits $E_0, \dots, E_t$ when $D$ is the constant $1$
as a proof of $C$ in a {\em rewrite proof system}: such proofs are p-time 
recognizable by Lemma \ref{ptime} as required by the Cook-Reckhow definition of 
proof systems \cite{CooRec}.
It follows that unless $\mbox{NP} = \mbox{coNP}$, given $C,D$ you cannot construct 
circuits $E_i$ satisfying all
the requirements (even using any other definition of $A \approx B$ that is p-time and 
that implies
the logical validity of $A \equiv B$) without some additional assumption (as
is the existence of a CF-proof of $C \equiv D$ with a bounded number of steps in the theorem).

\bigskip
\noindent
{\large {\bf Acknowledgments:}} I am indebted to Agni Datta (U. of Edinburgh) and
Archisman Dutta (Aarhus U.) for pointing out a mistake in an earlier version of the paper.

\end{document}